\documentclass[twoside,reqno]{amsart}

\usepackage{amsmath}
\usepackage{amssymb,latexsym,amstext,amsthm}
\usepackage{verbatim} 
\usepackage{eucal} 

\theoremstyle{plain}
\newtheorem{thm}[equation]{Theorem}
\newtheorem{pro}[equation]{Proposition}
\newtheorem{cor}[equation]{Corollary}
\newtheorem{lem}[equation]{Lemma}

\theoremstyle{definition}
\newtheorem{exa}[equation]{Example}

\newtheorem{DEF}[equation]{Definition}
\newtheorem{rem}[equation]{Remark}

\def\0b{\bar{0}}

\def\epb{\bar{\ep}}

\def\andd{\quad\hbox{and}\quad}
\def\ind{\hbox{ind}}

\def\v{{\mathcal V}}

\def\vd{\dot{\mathcal V}}
\def\vz{{\mathcal V}^{0}}
\def\vt{\tilde{\mathcal V}}

\def\fm{(\cdot,\cdot)}
\def\a{\alpha}

\def\w{{\mathcal W}}

\def\sub{\subseteq}
\def\rd{\dot{R}}

\def\lam{\lambda}
\def\Lam{\Lambda}

\def\1k{\frac{1}{k}}

\def\la{\langle}
\def\ra{\rangle}

\def\d{\delta}

\def\b{\beta}

\def\sg{\sigma}

\def\bbbz{{\mathbb Z}}

\def\ep{\epsilon}

\def\da{\dot{\a}}

\def\proof{\noindent{\bf Proof. }}

\def\rank{\hbox{rank}}

\def\a{\alpha}

\def\andd{\quad\hbox{and}\quad}

\def\b{\beta}

\def\d{\delta}

\def\e{\epsilon}

\def\ep{\epsilon}
\def\ind{\hbox{ind}}
\def\fm{(\cdot,\cdot)}

\def\sub{\subseteq}
\def\rd{\dot{R}}

\def\lam{\lambda}
\def\Lam{\Lambda}

\def\1k{\frac{1}{k}}

\def\la{\langle}
\def\ra{\rangle}

\def\qed{\hfill$\Box$\vspace{5mm}}
\def\sg{\sigma}

\def\rtimes{R^{\times}}

\def\vd{\dot{\mathcal V}}

\def\vz{{\mathcal V}^{0}}

\def\vt{\tilde{\mathcal V}}

\def\w{{\mathcal W}}

\def\bbbz{{\mathbb Z}}

\def\s{S}

\def\e{E}
\def\r{R}

\def\vb{\overline{\v}}

\def\supp{\hbox{supp}}

\def\FO{\hbox{FO}}
\def\O{\hbox{O}}

\def\da{\dot{\a}}
\def\hw{\hat{w}}
\def\th{\hat{t}}
\def\hz{\hat{z}}
\def\bi{{\bf I}}
\def\bii{{\bf II}}
\def\hu{\hat{u}}
\def\esupp{\hbox{Esupp}}
\def\mb{\tilde{m}}

\def\epb{\tilde{\ep}}

\begin{document}
\setcounter{page}{1}

\title{Presentation by conjugation for $A_1$-type\\ extended affine Weyl groups}

\author{Saeid Azam, Valiollah Shahsanaei}
\address
{Department of Mathematics\\ University of Isfahan\\Isfahan, Iran,
P.O.Box 81745-163} \email{azam@sci.ui.ac.ir, saeidazam@yahoo.com}
\thanks{The authors would like to thank the Center of Excellence for
Mathematics, University of Isfahan.}

\begin{abstract}
There is a well-known presentation for finite and affine Weyl
groups called the {\it presentation by conjugation}. Recently, it
has been proved that this presentation holds for certain
sub-classes of extended affine Weyl groups, the Weyl groups of
extended affine root systems. In particular, it is shown that if
nullity is $\leq 2$, an $A_1$-type extended affine Weyl group has
the presentation by conjugation. We set up a general framework for
the study of simply laced extended affine Weyl groups. As a
result, we obtain certain necessary and sufficient conditions for
an $A_1$-type extended affine Weyl group of arbitrary nullity to
have the presentation by conjugation. This gives an affirmative
answer to a conjecture that there are extended affine Weyl groups
which are not presented by ``presentation by conjugation".
\end{abstract}
\maketitle

\markboth{PRESENTATION BY CONJUGATION}{S. AZAM, V. SHAHSANAEI}

\medskip
\setcounter{section}{-1}
\section{\bf INTRODUCTION}\label{introduction}
There are several well-known presentations for finite and affine
Weyl groups (\cite{St}, \cite{H}, \cite{MP}) including the Coexter
presentation and the so called the {\it presentation by
conjugation} \cite{MP}. The elements of the Coexter presentation
can be read from the corresponding (type dependent) Cartan matrix
while the elements of the presentation by conjugation are given
uniformly for all types (see Definition \ref{DEF}). Our main focus
in this work is on the presentation by conjugation for $A_1$-type
{\it extended affine Weyl groups}.

Extended affine Weyl groups (Definition \ref{eawg}) are the Weyl
groups of extended affine root systems (Definition \ref{ears}). We
use the abbreviations EARS's and EAWG's for extended affine root
systems and extended affine Weyl groups, respectively. Finite and
affine root systems are nullity zero and nullity one extended
affine root systems, respectively. In 1995, Y. Krylyuk \cite{K}
showed that simply laced EAWG's of rank $>1$ have the presentation
by conjugation. He used this to establish certain classical
results for extended affine Lie algebras. In 2000, S. Azam
\cite{A3} proved that the presentation by conjugation holds for a
large subclass of EAWG's, including reduced EAWG's of nullity
$\leq 2$.

In this work we obtain certain necessary and sufficient conditions
for arbitrary nullity EAWG's of type $A_1$ to have the
presentation by conjugation. This has several interesting
implications which are recorded in Section \ref{type-1}. In
Section \ref{semilattices}, we recall the definition of a {\it
semilattice} from \cite{AABGP} and introduce the notions of {\it
support}, {\it essential support} and {\it integral collections}
for semilattices. These notions naturally arise in the study of
EAWG's and their presentations. In Sections \ref{EAWG} and
\ref{results}, we record some important results mainly from
\cite{AS} which are crucial for the up coming sections, and
establish a few preliminary results regarding EAWG's.

In Section \ref{by-conj}, we assign a presented group $\hat{\w}$
to the Weyl group $\w$ of an extended affine root system $R$
(Definition \ref{DEF}) and study its structure in details.
Following \cite{K}, we say the Weyl group $\w$ has the
presentation by conjugation if $\w\cong\hat{\w}$. Since our
setting has set up for the study of general simply laced EARS's,
the results of this section gives in part a simple and new proof
of \cite[Theorem III.1.14]{K}, namely simply laced EAWG's of rank
$>1$ have the presentation by conjugation. In Section
\ref{type-1}, the main section, we compute the kernel of a natural
epimorphism $\psi:\hat{\w}\longrightarrow\w$ and show that it is
isomorphic to the direct sum of a finite number of copies of the
cyclic group of order $2$ (Proposition \ref{redu-3}). This then
allows us to prove that $\w$ has the presentation by conjugation
if and only if the center of $\hat{\w}$ is free abelian group, if
and only if a certain integral condition (see Section
\ref{semilattices} for definition) holds for the semilattice
involved in the structure of $R$, if and only if a minimal
condition holds on a particular set of generators for the Weyl
group (Theorem \ref{reduced}). This reduces the problem of
deciding which EAWG has the presentation by conjugation to a
finite problem which can in principle be solved for any $\nu\geq
0$. Using this criterion it is then easy to show that if $\nu\leq
3$, the Weyl groups of all EARS's of type $A_1$ have the
presentation by conjugation except only one root system, namely
the one of index $7$ which happens in nullity $3$ (Corollary
\ref{main}). This in particular gives an affirmative answer to a
conjecture due to S. Azam \cite[Remark 3.14]{A3} that certain
EAWG's of type $A_1$ may not have the presentation by conjugation
(see Remark \ref{rem1}).

For the study of EARS's and EAWG's we refer the reader to
\cite{MS}, \cite{Sa}, \cite{AABGP}, \cite{A1}, \cite{A2},
\cite{A3}, \cite{A4}, \cite{SaT}, \cite{T} and \cite{AS}. The
authors would like to thank A. Abdollahi for a fruitful discussion
which led to the proof of crucial Lemma \ref{minimal}. The authors
hope that the results appearing in this work prepare the ground
for the study of the presentation by conjugation for non-simply
laced extended affine Weyl groups.

\section{\bf SEMILATTICES}
\label{Semilattices} \setcounter{equation}{0}\label{semilattices}
In this section, we briefly recall the definition of a semilattice
from \cite{AABGP} and record certain properties of semilattices
which will be important for our work. We also introduce a notion
of {\it integral collections} for semilattices which turn out to
be crucial in the study of Weyl groups of extended affine root
systems. For the details on semilattices the reader is referred to
\cite[Chapter II.\S 1]{AABGP}. In this section, we fix several
sets which will be used throughout the paper.

\begin{DEF}\label{semi1} A {\it semilattice} is a subset $S$ of a
finite dimensional real vector space $\v^0$ such that $0\in S$,
$S\pm 2S\sub S$, $S$ spans $\v^0$ and $S$ is discrete in $\v^0$.
The {\it rank} of $S$ is defined to be the dimension $\nu$ of
$\v^0$. Note that the replacement of $\s\pm2\s\subseteq \s$ by
$\s\pm\s\subseteq \s$ in the definition gives one of the
equivalent definitions for a {\it lattice} in $\vz$.
\end{DEF}

Let $S$ be a semilattice in $\v^0$. The $\bbbz$-span $\Lam$ of $S$
in $\v^0$ is a lattice in $\v^0$ with a basis $B$ consisting of
elements of $S$. Namely
\begin{equation}\label{basis-1}
B=\{\sg_1,\ldots,\sg_\nu\}\sub S\quad\hbox{with}\quad
\Lam=\sum_{i=1}^{\nu}\bbbz\sg_i. \end{equation}
 We fix this basis
$B$. For a set $J\sub\{1,\ldots,\nu\}$ we put
$$\tau_{_J}:=\sum_{j\in J}\sg_j.$$
(If $J=\emptyset$ we have by convention $\sum_{j\in J}\sg_j=0$).
By [AS, \S 1], there is a unique subset, denoted $\supp(S)$,
consisting of subsets of $\{1,\ldots,\nu\}$ such that
\begin{equation}\label{unique}
\emptyset\in\supp(S)\andd
S=\bigcup_{J\in\supp(S)}(\tau_{_J}+2\Lam).
\end{equation}
Note that since $B\sub S$, we have $\{r\}\in\supp(S)$ for all
$r\in J_\nu$.  The collection $\supp(S)$ is called the {\it
supporting class of} $S$ (with respect to $B$). We also introduce
a notion of {\it essential support} for $S$, denoted $\esupp(S)$,
namely
$$\esupp(S)=\{J\in\supp(S):\;|J|\geq 3\}.
$$
Clearly we have
$$|\esupp(S)|\leq 2^{\nu}-1-\nu-\left(
\begin{array}{c}
  \nu \\
  2 \\
\end{array}\right).
$$
Following \cite{A4}, we call the integer $\ind(S):=|\supp(S)|-1$,
{\it index} of $S$.

We assign to each semilattice $S$ certain collections of integers,
called {\it integral collections}. This notion will play a crucial
role in the sequel.

 Let $S$ be a semilattice of rank $\nu$.
 For $1\leq r< s\leq\nu$ and $J\in\supp(S)$ we set
\begin{eqnarray*}
 \d(r,s)=\left\{\begin{array}{ll}
 1,& \hbox{if }\{r,s\}\in\supp(S),\vspace{2mm}\\
 2, &\hbox{if }\{r,s\}\not\in\supp(S),\vspace{2mm}\\
  \end{array}\right. \andd\d(_J,r,s)=\left\{\begin{array}{ll}
 1,& \hbox{if }\{r,s\}\subsetneq J,\vspace{2mm}\\
 0, &\hbox{otherwise}.
 \end{array}\right.
   \end{eqnarray*}
We call a collection $\epb=\{\ep_{_J}\}_{J\in\esupp(S)}$,
$\ep_{_J}\in\{0,1\}$, an {\it integral collection} for $S$ if
$$
\frac{1}{\d(r,s)}\sum_{J\in\esupp(S)}
\d(_J,r,s)\epsilon_{_J}\in\bbbz,\quad\hbox{for all}\quad 1\leq
r<s\leq\nu.
$$
If $\esupp(S)=\emptyset$, we interpret $\epb$ to be the zero
collection and the sum on an empty set to be zero. If $\ep_{_J}=0$
for all $J\in\esupp(S)$, we call $\epb=\{0\}_{J\in\esupp(S)}$ the
{\it trivial} collection. Clearly the trivial collection is an
integral collection and there are at most $2^{|\esupp(S)|}$
integral collections for $S$. Any integral collection different
from the trivial collection is called {\it non-trivial}.

\begin{exa}\label{exa1}
(i) If $\esupp(S)=\emptyset$ or $\ep_{_J}=0$ for all
$J\in\esupp(S)$, then it is clear from definition that the only
integral collection for $S$ is the trivial collection.

(ii) Let $S$ be such that $\{r,s\}\in\supp(S)$ for all $1\leq
r<s\leq\nu$. Then $\d(r,s)=1$ for all such $r<s$ and so any
collection $\{\ep_{_J}\}_{J\in\esupp(S)}$, $\ep_{_J}\in\{0,1\}$ is
an integral collection. Therefore there are $2^{|\esupp(S)|}$
integral collections for $S$. In particular, if $S$ is a lattice
of rank $\nu\geq 3$, then there is at least one non-trivial
integral collection for $S$, as $|\esupp(S)|\geq 1$.


(iii) Suppose $\esupp(S)$ contains a set $J'$ such that
$\{r,s\}\in\supp(S)$ for all $r,s\in J'$. Set
$$
\ep_{_J}=\left\{\begin{array}{ll}
  1, &\hbox{if }J=J',\vspace{2mm} \\
  0 & \hbox{otherwise}.
\end{array}\right.$$
Now for $1\leq r<s\leq\nu$ we have $\d(J',r,s)=0$ if
$\{r,s\}\not\sub J'$. It follows that $\{\ep_{_J}\}$ is a
non-trivial integral collection for $S$. In particular, if $\nu>3$
and $\ind(S)=2^\nu-2$, then there is at least one non-trivial
integral collection for $S$, as in this case there is always a set
$J'\in\esupp(S)$ with $\{r,s\}\in\supp(S)$ for all $r,s\in J'$.

(iv) Let $S$ and $S'$ be semilattices with the same $\bbbz$-span
$\Lam$ and with $S\sub S'$. Then $\supp(S)\sub\supp(S')$. Now if
$\epb$ is a non-trivial integral collection for $S$ then the
collection $\epb'=\{\ep'_{_J}\}_{J\in\esupp(S')}$ defined by $$
\ep'_{_J}=\left\{
\begin{array}{ll}
\ep_{_J}&\hbox{if }J\in\esupp(S)\\
0&\hbox{otherwise}. \end{array}\right. $$ is a non-trivial
integral collection for $S'$.
\end{exa}

\begin{lem}\label{reduced-1}
 If $\ind(S)-\nu\leq3$,
 then the only integral collection for $S$ is
the trivial collection.
 \end{lem}

\proof By Example \ref{exa1}(i), it is enough to show that if
$\{\ep_{_J}\}$ is an integral collection for $S$ with
$\esupp(S)\not=\emptyset$ then $\ep_{_J}=0$ for all
$J\in\esupp(S)$. Now let $J\in\esupp(S)$. Since
$\{\emptyset,\{1\},\ldots,\{\nu\}\}\sub\supp(S)$ and
$\ind(S)\leq\nu+3$, it is clear that $\esupp(S)$ has at most three
elements. Thus we may choose $r<s$ such that $\{r,s\}\sub J$,
$\{r,s\}\not\in\supp(S)$ and $\{r,s\}\not\sub J'$ for all
$J'\in\esupp(S)$ with $J'\not=J$. Then from the definition of an
integral collection we get $\ep_{_J}/2\in\bbbz.$ This forces
$\ep_{_J}=0$.\qed

\section{\bf EXTENDED AFFINE WEYL GROUPS}
\label{EAWG} \setcounter{equation}{0} In this section, we briefly
recall the definition and some basic facts about extended affine
root systems (EARS for short) an their corresponding Weyl groups.
For a more detailed exposition, we refer the reader to [AABGP] and
[AS], in particular we will use the notation and concepts
introduced there without further explanations.

All the groups we consider in this work will be subgroups of the
orthogonal group $\O(\vt,I)$ where $\vt$ is a finite dimensional
vector space equipped with a non-degenerate symmetric bilinear
form $I=\fm$. For a group $G$, $Z(G)$ denotes the center of $G$,
and for $x,y\in G$ we denote the commutator $x^{-1}y^{-1}xy$ by
$[x,y]$. An element $\a\in\vt$ is called {\it non-isotropic} ({\it
isotropic}) if $(\a,\a)\not=0$ ($(\a,\a)=0$). We denote the set of
non-isotropic elements of a subset $T$ with $T^\times$ and the set
of isotropic elements of $T$ with $T^0$.

\begin{DEF}\label{ears}
{\rm A subset $R$ of a non-trivial finite dimensional real vector
space $\v$, equipped with a positive semi-definite symmetric
bilinear form $(\cdot,\cdot)$, is called an {\it extended affine
root system} belonging to $(\v,(.,.))$ (EARS for short) if $\r$
satisfies the following 8 axioms:
\begin{itemize}

\item  R1) $ 0\in \r$,

\item R2) $ -\r=\r$,

\item R3) $\r$ spans $\v$,

\item R4) $\a \in \r^{\times} \Longrightarrow 2\a \notin \r$,

\item R5) $\r$ is discrete in $\v$,

\item R6) For $\a \in \r^{\times}$  and $\beta\in \r$, there exist
non-negative integers $d,u$ such that $\b+n\a\in R$, $n\in\bbbz$,
if and only if $-d\leq n\leq u$, moreover $2(\b,\a)/(\a,\a)=d-u$.

\item R7) If $\r=\r_{1}\cup \r_{2}$, where $(\r_{1},\r_{2})=0$,
then either $\r_{1}=\emptyset$ or $\r_{2}=\emptyset$,

\item  R8)   For any $\sg \in \r^{0}$, there exists $\a\in
\r^{\times}$ such that $\a+\sg\in R$.
\end{itemize}}
\end{DEF}

The dimension $\nu$ of the radical $\v^0$ of the form is called
the {\it nullity} of $R$, and the dimension $\ell$ of
$\vb:=\v/\v^0$ is called the {\it rank} of $R$. It turns out that
the image of $R$ in $\vb$ is a finite root system whose type is
called the {\it type} of $R$. Corresponding to the integers $\ell$
and $\nu$, we set
$$J_\ell=\{1,\ldots,\ell\}\andd J_\nu=\{1,\ldots,\nu\}.
$$

Let $\r$ be a $\nu$-EARS  belonging to $(\v,(.,.))$. It follows
that the form restricted to $\vb$ is positive definite and that
$\bar{\r}$, the image of $R$ in $\vb$, is an irreducible finite
root system (including zero) in $\vb$ ([AABGP, II.2.9]). The {\it
type} of $R$ is defined to be the type of $\bar{R}$.

In this work {\it we always assume that $R$ is an EARS of simply
laced type}, that is it has one of the types $X_\ell=A_\ell$,
$D_\ell$, $E_6$, $E_7$ or $E_8$. We fix a complement $\vd$ of
$\v^0$ in $\v$ such that
$$
\rd:=\{\da\in\vd\mid\da+\sg\in R\hbox{ for some }\sg\in\v^0\}
$$
is a finite root system in $\vd$, isometrically isomorphic to
$\bar{R}$, and that
\begin{equation}\label{AABGP}
R=R(X_\ell,S)=(S+S)\cup(\rd+S)
\end{equation} where $S$ is
a semilattice in $\v^0$. Recall that we have fixed a basis $B$ of
$\v^0$ satisfying (\ref{basis-1}). Here $X_\ell$ denotes the type
of $\rd$. Throughout this work {\it we fix a fundamental basis}
$$\dot{\Pi}=\{\a_1,\ldots,\a_\ell\}$$ of $\rd$. It is known that if
$\ell>1$, then $S$ is a lattice in $\v^0$. To introduce the
extended affine Weyl group of $R$, we need to consider the so
called {\it a hyperbolic extension} $\vt$ of $\v$ as follows. We
fix a basis $\{\lam_1,\ldots,\lam_\nu\}$ of $(\v^0)^*$ and we set
$\vt:=\v\oplus (\vz)^*,$ where  $(\vz)^*$ is the dual space of
$\vz$. Now we extend the from $\fm$ on $\v$ to a non-degenerate
form $I$, denoted again by $\fm$, on $\vt$ as follows:
$$
\begin{array}{ll}
 \bullet& (.,.)_{|_{\v\times\v}}:=(.,.),\vspace{1mm}\\
\bullet& (\vd,(\vz)^*)=((\vz)^*,(\vz)^*):=0,\vspace{1mm}\\
\bullet& (\sg_i,\lam_j):=\delta_{i,j},\quad i,j\in J_\nu.
\end{array}
$$

\begin{DEF}\label{eawg} The {\it extended affine Weyl group} $\w$ of
$R$ is the subgroup of the orthogonal group $\O(\vt,I)$ generated
by reflections $w_\a$, $\a\in\rtimes$, defined by
$w_\a(u)=u-(u,\a)\a$, $u\in\vt$. Clearly
\begin{equation}\label{fixed}
\w\sub\FO(\vt,I):=\{w\in \O(\vt,I)\mid w(\d)=\d\hbox{ for all
}\d\in\v^0\}.
\end{equation}
The subgroup
$${\mathcal H}:=\la
w_{\a+\sg}w_\a~~\mid\a\in \r^\times,\;\sg\in\vz, \a+\sg\in\r\ra$$
of $\w$ is called the {\it Heisenberg-like group} of $R$ and plays
a crucial role in the study of structure and presentations of
$\w$. Using the fact that $\rd\sub R$, we identify the finite Weyl
group $\dot{\w}$ of $\rd$ with a subgroup of $\w$. Then for any
$\da\in\rd^\times$, using (\ref{fixed}), (\ref{AABGP}),
(\ref{unique}) and the fact that $\rd^\times=\dot{\w}\da$ we have
 \begin{eqnarray}\label{union}
 \r^\times= \bigcup_{J\in\supp(S)}(\dot{\w}(\da)+\tau_{_J}+2\Lam)=
 \dot{\w}\big(\bigcup_{J\in\supp(S)}(\da+\tau_{_J}+2\Lam)\big).
 \end{eqnarray}
 We also note that the following important conjugation relation
 holds in $\w$:
 \begin{equation}\label{normal}
 ww_\a w^{-1}=w_{w\a}\qquad (\a\in\rtimes,\;w\in\w).
 \end{equation}
\end{DEF}

In [AS, \S 3] we have studied the structure of $\w$ and have
obtained a particular finite set of generators for $\w$ and its
center $Z(\w)$, using linear maps defined as follows. For
$\a\in\v$ and $\sg\in\vz$, we define $T^{\sg}_{\a}\in\hbox{\e
nd}(\vt)$ by
\begin{equation}\label{xi}
T^{\sg}_{\a}(u):=u-(\sg,u)\a+(\a,u)\sg
-\frac{(\a,\a)}{2}(\sg,u)\sg\quad(u\in\vt).
\end{equation}
One can check easily that for $\a\in\v^\times$,\; $\b\in\v$\;
$\sg,\d\in\vz$ and $w\in\mbox{O}(\vt,I)$,
\begin{equation}\label{translation}
\begin{array}{c}
   T_{\a}^{\sg}=w_{\a+\sg}w_{\a},\quad [T^{\sg}_\a,T^{\d}_\b]=T_\sg^{(\a,\b)\d}
     \andd    wT^{\d}_\b w^{-1}=T^{\d}_{w(\b)}.
\end{array}
 \end{equation}
For $\a\in\v$ and $\sg,\d\in\v^0$, we have from
(\ref{translation}) that
\begin{equation}\label{neww}
 T^\sg_\a\in
\hbox{FO}(\vt,I)\andd T^\sg_\d\in Z\big(\hbox{FO}(\vt,I)\big).
\end{equation}
For  $r,s\in J_\nu$, we set
\begin{equation}\label{def}
 c_{r,s}=T_{\sg_r}^{\sg_s}\andd C:=\la
c_{r,s}:1\leq r<s\leq\nu\ra.
\end{equation}
 Then by
(\ref{neww}) for all $r,r',s,s'\in
   J_{\nu}$ and $\a\in\v$, we have
\begin{eqnarray}\label{ex-3}
  c_{r,s}(\a)=\a\andd c_{r,r}=c_{r,s}c_{s,r}=[c_{r,s},c_{r',s'}]=1.
 \end{eqnarray}
By \cite[Lemma 2.16(ii)]{AS},
\begin{equation}\label{free-abelian}
C\hbox{ is a free abelian group of rank }{\nu(\nu-1)}/{2}.
\end{equation}
Also for $(i,r)\in J_\ell\times J_\nu$, we set
\begin{equation}\label{newdef}
 t_{i,r}:=T_{\a_i}^{\sg_r}=w_{\a_i+\sg_r}w_{\a_i}.
\end{equation}
From (\ref{translation}) for all $r,s\in J_{\nu}$, $i,j\in J_\ell$
and $\a\in\v$, one can see easily that
\begin{eqnarray}\label{newex-3}
[t_{i,r},t_{j,s}]=c_{r,s}^{(\a_i,\a_j)}\andd w_{\a_j}t_{i,r}
w_{\a_j}=t_{i,r}t_{j,r}^{-(\a_i,\a_j)}.
 \end{eqnarray}
To describe the center $Z(\w)$ of $\w$,  we set
$$Z:=\la z_{_J}~|~J\subseteq J_\nu\ra\leqq Z\big(\hbox{FO}(\vt,I)\big),$$
where
\begin{equation}\label{def-z}
z_{_J}:=\left\{\begin{array}{ll} \prod_{\{r,s\in J\mid r<s\}}
  c_{_{r,s}} & \hbox{if }J\in\supp(S)\vspace{2mm}\\
   c_{_{r,s}}^{2} &  \hbox{if }J=\{r,s\}\not\in\supp(S),\; r< s\vspace{2mm}\\
   1 & \hbox{otherwise}.
\end{array}\right.
\end{equation}
(Here we interpret the product on an empty index set to be $1$.)
We note from the definitions of $z_{_J}$ and $\d(r,s)$ that
\begin{equation}\label{final1}
z_{_{\{r,s\}}}=c_{r,s}^{\d(r,s)}\andd
c_{r,s}^2=z_{_{\{r,s\}}}^{3-\d(r,s)}.
\end{equation}

\begin{lem}\label{indepen}
  $\la z_{_{\{r,s\}}}~|~1\leq r<s\leq\nu\ra$ is
 a free abelian group of rank $\nu(\nu-1)/2$.
 \end{lem}
 \proof From (\ref{def}) and (\ref{def-z}), we see that
 the group in the statement is squeezed
 between two groups $\la c_{r,s}^2:1\leq
r<s\leq\nu\ra$ and $C$. Since $C$ is free abelian on generators
$c_{r,s}$, $1\leq r<s\leq\nu$, the result follows.\qed

Using \cite[Propositions 3.5, 3.13, 3.14, 3.16]{AS} we have the
following two propositions.
\begin{pro}\label{conj-6}
(i) $\w=\dot{\w}\ltimes\mathcal H$.

(ii) $\w=\la w_{\a_i}, t_{i,r},z_{_J}\mid i\in J_\ell,\;r\in
  J_\nu,\;J\sub J_\nu\ra.$

(iii) $\mathcal H=\la  t_{i,r},z_{_J}\mid i\in J_\ell,\;r\in
  J_\nu,\;J\sub J_\nu\ra.$

(iv) $Z(\w)=Z(\mathcal H)=Z.$

(v) $\mathcal H$ is a torsion free, two step nilpotent group.

(vi) $Z$ is
 a free abelian group of rank $\nu(\nu-1)/2$.
\end{pro}
\begin{pro}\label{unique-6}
 Each element $w\in \w$ has a unique expression of the form
\begin{equation}\label{exp} w=z\dot{w}
\prod_{i=1}^{\ell}\prod_{r=1}^{\nu}t_{i,r}^{m_{i,r}}
  \quad(\dot{w}\in\dot{\w},\;z\in Z(\w),\; m_{i,r}\in\bbbz).
\end{equation}
\end{pro}
\medskip
\section{\bf PRELIMINARY RESULTS}\label{results}\setcounter{equation}{0}
We keep all the notations as in Section \ref{EAWG}. In particular
$R$ is an EARS of simply laced type and $\w$ is its extended
affine Weyl group. We have fixed a basis $B$ of $\v^0$ and a
finite root system $\rd$ such that $R$ is of the form
(\ref{AABGP}), where the semilattice $S$ is given by
(\ref{unique}). Also recall from Section \ref{EAWG} that we have
fixed a basis $\dot{\Pi}=\{\a_1,\ldots,\a_\ell\}$ of $\rd$.

For a subset  $J=\{i_1,\ldots, i_n\}$ of $J_\nu$ with $
i_1<i_2<\cdots <i_n$ and a group $G$ we make the convention
$$\prod_{i\in J} a_i=a_{i_1}a_{i_2}\cdots a_{i_n}\qquad\qquad(a_i\in G). $$
Using this notation, (\ref{xi}) and (\ref{newdef}), we have
\begin{equation}\label{prod1}
\prod_{r\in J}t^{m_r}_{i,r}(\a_j)=\a_j+\sum_{r\in
J}(\a_i,\a_j)m_r\sg_r \qquad(\a_i,\a_j\in\dot{\Pi},\;m_r\in\bbbz).
\end{equation}
\begin{lem}\label{conj-sim}
Let $J\sub J_\nu$ and $i\in J_\ell$.

(i) If $J\in\mbox{supp}(S)$ then
$$z_{_J}=w_{\a_i}w_{\a_i+\tau_{_J}}\prod_{r\in
    J}w_{\a_i+\sg_r}w_{\a_i}.$$

(ii) If $J=\{r,s\}\not\in\mbox{supp}(S)$, $r<s$, then
$$z_{_J}=[t_{i,r}, t_{i,s}]=w_{\a_i}w_{\a_i+\sg_r}
w_{\a_i}w_{\a_i+\sg_s}w_{\a_i+\sg_r}w_{\a_i}w_{\a_i+\sg_s}w_{\a_i}.$$
\end{lem}

\proof (i) From (\ref{translation}),
  \cite[Lemma 2.8]{AS} and the way $z_{_J}$ is defined  we have
\begin{eqnarray*}
w_{\a_i+\tau_{_J}}w_{\a_i}=T_{\a_i}^{\tau_{_J}}=T_{\a_i}^{\sum_{r\in
J}\sg_r}=\prod_{r\in J}t_{i,r}\prod_{\{r,s\in J\mid
r<s\}}c_{s,r}=\prod_{r\in J}t_{i,r}z_{_J}^{-1}.
\end{eqnarray*}
Then \begin{eqnarray*} z_{_J}=w_{\a_i}w_{\a_i+\tau_{_J
}}\prod_{r\in J}t_{i,r}=w_{\a_i}w_{\a_i+\tau_{_J}}\prod_{r\in
    J}w_{\a_i+\sg_r}w_{\a_i}.
\end{eqnarray*}

(ii) It is an immediate consequence of the way $z_{_{\{r,s\}}}$ is
defined and (\ref{newex-3}). \qed

We now consider a subset $\Pi$ of $\r^\times$ which in some
aspects behave similar to a basis for a finite root system. The
set $\Pi$ is introduced in \cite[\S 4]{A4} for reduced extended
affine root systems and its properties is studied. Set
\begin{equation}\label{pi}
\Pi=\left\{\begin{array}{ll}
\big\{\a_1+\tau_{_J}\mid J\in\supp(S)\}, & \hbox{if }\ell=1, \vspace{2mm} \\
   \dot{\Pi}\cup\big\{\a_1+\sg_r\mid r\in J_\nu\big\}, & \hbox{if }\ell>1,
 \end{array}\right.
\end{equation}
 and
 $$\w_{\Pi}=\la w_\a\mid \a\in\Pi\ra.$$
 It is clear that $\dot{\w}\sub\w_{\Pi}$.

 \begin{lem}\label{belong}
 $t_{i,r}\in\w_{\Pi}$, $(i,r)\in J_\ell\times J_\nu$.
 \end{lem}
 \proof Let $w\in\dot{\w}$ be such that
$\a_i=w(\a_1)$. Since $\dot{\w}\sub\w_{\Pi}$ and
$\{r\}\in\supp(S)$ for all $r\in J_\nu$, we have from
(\ref{normal})
$$t_{i,r}=w_{\a_i+\sg_r}w_{\a_i}=w_{w(\a_1)+\sg_r}w_{w(\a_1)}=
ww_{\a_1+\sg_r}w_{\a_1}w^{-1}\in\w_{\Pi}.$$ \qed

The following result is proved in \cite[Proposition 4.26]{A4}. We
present a new proof for this result here.

\begin{lem}\label{base-1}
 $\w_{\Pi}\Pi=\r^\times$ and $\w=\w_{\Pi}$.
 \end{lem}
 \proof It is clear that $\w_{\Pi}\Pi\subseteq
 \r^\times$. So to prove the first equality, we must show that  $\r^\times \subseteq\w_{\Pi}\Pi$.
 Using (\ref{union}) and the fact that $\dot{\w}\sub\w_{\Pi}$, it is enough to show that for  $J\in\supp(S)$ and
$m_r\in\bbbz$, $r\in J_\nu$,
$$\a:=\a_1+\tau_{_J}+\sum_{r=1}^\nu2m_r\sg_r\in\w_{\Pi}\Pi.$$
By (\ref{prod1}), (\ref{fixed}) and the fact that
$\tau_{_J}=\sum_{i\in J}\sg_i$, we have
$$
\a=\prod_{r=1}^{\nu} t_{1,r}^{m_r}(\a_1+\tau_{_J}),
$$
which if $\ell=1$ is an element of $\w_\Pi\Pi$, by Lemma
\ref{belong}. If $\ell>1$, then there exists $\a_j\in\dot{\Pi}$
such that $(\a_1,\a_j)=-1$. Thus
$$\a=\prod_{r=1}^{\nu} t_{1,r}^{m_r}(\a_1+\tau_{_J})
=\prod_{r=1}^{\nu}t^{m_r}_{1,r}\prod_{s\in
J}t_{j,s}^{-1}(\a_1+\sg_s),
$$
which is again an element of $\w_\Pi\Pi$, by Lemma \ref{belong}.
Thus the first equality in the statement holds. Using
(\ref{normal}) and the first equality, one can see that  the
second equality holds.\qed

 \section{\bf PRESENTATION
 BY CONJUGATION}\label{by-conj}\setcounter{equation}{0}

We keep all the notations and assumptions  as in the previous
sections. In particular, $R$ is a nullity $\nu$ simply laced
extended affine root system of the form (\ref{AABGP}), where the
semilattice $S$ is given by (\ref{unique})   and $\w$ is its
extended affine Weyl group.

\begin{DEF}\label{DEF} Let $\hat{\w}$ be the group defined by generators
$\hw_{\a}$, $\a\in\r^\times$ and relations:
$$\begin{array}{l}
  (\bi) \hspace{2mm}\hw_\a^2=1,\qquad \a\in\r^\times\vspace{2mm} \\
  (\bii)\hspace{2mm}\hw_\a\hw_\b\hw_\a=\hw_{w_{\a}(\b)},\qquad
   \a,\b\in\r^\times. \\
\end{array}$$
Following \cite{K}, we say that the extended affine Weyl group
$\w$ of $R$ has {\it the presentation by conjugation} if
$\w\cong\hat{\w}$. \end{DEF}

By (\ref{normal}) and the fact that $w_{\a}^2=1$ for
$\a\in\r^\times$, the relations of the forms (\bi) and (\bii) are
satisfied in $\w$ (replacing $\hat{w}_{\a}$ with $w_{\a}$ for any
$\a\in\r^\times$). So  the assignment $\hw_{\a}\longmapsto w_{\a}$
induces a unique epimorphism
 \begin{equation}\label{main}
 \psi:\hat{\w}\longrightarrow \w.
 \end{equation}
Since the finite Weyl group $\dot{\w}$ has the presentation by
conjugation \cite{St}, the restriction of $\psi$ to
$\hat{\dot{\w}}:=\la\hw_{\a_i}\mid i\in J_\ell\ra$ induces the
isomorphism
 \begin{equation}\label{finite-case}
 \hat{\dot{\w}}\cong^{^{\hspace{-2mm}\psi}}\dot{\w}.
 \end{equation}
  One can easily deduce from relations (\bi) and (\bii) that
\begin{equation}\label{main1}
\hw\hw_\a\hw^{-1}=\hw_{\psi(\hw)(\a)}\qquad\qquad(\hw\in\hat{\w},\;\a\in\rtimes).
\end{equation}

 For any $(i,r)\in J_\ell\times J_\nu$
  and $J\subseteq J_\nu$ we set
$$\th_{i,r}:=\hw_{\a_i+\sg_r}\hw_{\a_i},$$
and
\begin{equation}\label{test-1}
\hz_{_J}=\left\{\begin{array}{ll}
  \hw_{\a_1}\hw_{\a_1+\tau_{_J}}\prod_{r\in
    J}\th_{1,r},& \hbox{if}\;\;J\in\supp(S),\vspace{2mm} \\
  \big[\th_{1,r},\;\th_{1,s}], & \hbox{if}\;\;J=\{r,s\}\not\in\supp(S),\;r<s, \vspace{2mm} \\
  1&\hbox{otherwise.}
\end{array}\right.
\end{equation}
Since $\{r\}\in\supp(S)$ for $1\leq r\leq\nu$ and
$\a_1+\tau_{_J}\in\rtimes$ for all $J\in\supp(S)$, we have
\begin{equation}\label{test}
\la\hw_{\a_i},\;\th_{i,r},\;\hz_{_J}\mid(i,r)\in J_\ell\times
J_\nu,\;J\sub J_\nu\ra\sub\hat{\w}. \end{equation} Moreover,
because of (\bi), we get $\hz_{_{\{r\}}}=1$ for all $r\in J_\nu$.
Then by (\ref{newdef}) and Lemma \ref{conj-sim}, we have
\begin{equation}\label{psi}
\psi(\th_{i,r})=t_{i,r}\andd\psi(\hz_{_J})=z_{_J}.
\end{equation}

 \begin{lem}\label{generators-con-1}
(i) $\hat{\w}=\la\hw_\a\mid\a\in\Pi\ra.$

(ii) $\hat{\w}=\la\hw_{\a_i},\;\hat{t}_{i,r},\;\hz_{_J}\mid
(i,r)\in J_\ell\times J_\nu,\;J\sub J_\nu\ra.$

(iii) If  $\ell>1$, then $\hat{\w}=\la
\hw_{\a_i},\;\hat{t}_{i,r},\;\hz_{_{\{r,s\}}}\mid i\in
J_\ell,\;\;1\leq r<s\leq\nu\ra$.
\end{lem}

\proof (i) Let $\Pi$ be as in (\ref{pi}) and $\a\in\r^\times$.
From Lemma \ref{base-1} we have $\a=w(\b)$ for some $w\in\w_{\Pi}$
and $\b\in\Pi$. So $w=w_{\b_1}\cdots w_{\b_n}$, for some
$\b_i\in\Pi$. Then from (\bii) we have
$$\hw_\a=\hw_{w(\b)}=\hw_{\b_1}\cdots\hw_{\b_n}\hw_\b
\hw_{\b_n}\cdots\hw_{\b_1}\in\la\hw_\a\mid\a\in\Pi\ra.$$

(ii)-(iii) Considering  (\ref{test}) and part (i) we only need to
show that for any $\b\in\Pi$, $\hw_\b$ is in the right hand side
of the equality in the statement. Clearly this holds if $\b=\a_i$
for some $i\in J_\ell$. If $\b=\a_1+\tau_{J}$ for some
$J\in\supp(S)$, then from the way $\hz_{_J}$ is defined it is
clear that (ii) holds. Therefore to see (iii), it is enough to
show that $\hz_{_J}$ is in the right hand side of the equality in
(ii) for all $J\sub J_\nu$. But this clearly holds if
$J\not\in\supp(S)$. Now let $J\in\supp(S)$. Since $\ell>1$, there
exists $i\in J_\ell$ such that $(\a_1,\a_i)=-1$. Let  $\hat
w=\Pi_{r\in J}\th^{-1}_{i,r}$ and $w=\Pi_{r\in J}t^{-1}_{i,r}$. We
have from (\ref{main1}), (\ref{psi}) and (\ref{prod1}) that
\begin{eqnarray}\label{spi}
\hw_{\a_1+\tau_{_J}}=\hw_{w(\a_1)}=\hw_{\psi(\hat
w)(\a_1)}=\hat{w}\hat{w}_{\a_1}\hat{w}^{-1}=\prod_{r\in
J}\th^{-1}_{i,r}\hw_{\a_1}\big( \prod_{r\in
J}\th^{-1}_{i,r}\big)^{-1}.
\end{eqnarray} Therefore
$\hz_{_J}$ is in the right hand side and we are done.\qed

\begin{lem}\label{conj-cen-1}
(i) $[\hat{t}_{i,r},\hat{t}_{j,s}]\in Z(\hat{\w}),$ for
$(r,i),\;(j,s)\in J_\ell\times J_\nu$.

(ii) $\hz_{_J}\in Z(\hat{\w}),$ for $J\subseteq J_\nu$.
\end{lem}

\proof Let $\a\in\r^\times$.   From (\ref{main1}), (\ref{psi}),
(\ref{ex-3}), (\ref{newex-3}) and the definition of $z_{_J}$ we
have
$$[\hat{t}_{i,r},\hat{t}_{j,s}]\hw_\a[\hat{t}_{i,r},\hat{t}_{j,s}]^{-1}=
\hw_{[t_{i,r},t_{j,s}](\a)}=\hw_\a$$ and
$$\hz_{_J}\hw_\a\hz_{_J}^{-1}=\hw_{z_{_J}(\a)}=\hw_{\a}.$$
Thus $[\hat{t}_{i,r},\hat{t}_{j,s}]$ and  $\hz_{_J}$ commute with
all elements of $\hat{\w}$.\qed

\begin{lem}\label{well-define}
Let $i\in J_\ell$. Then

(i) $[\hat{t}_{1,r},\hat{t}_{1,s}]=[\hat{t}_{i,r},\hat{t}_{i,s}],$
for all $r,s\in J_\nu$,

(ii) $\hz_{_J}=\hw_{\a_i}\hw_{\a_i+\tau_{_J}}\prod _{r\in
    J}\th_{i,r}$ for all $J\in\supp(S)$. In particular,
if $\ell>1$, this holds for all $J\sub J_\nu$.
\end{lem}

\proof Let $w\in\dot{\w}$ be such that $\a_i=w(\a_1)$ and fix a
primage $\hw\in\hat{\w}$ of $w$, under $\psi$. Then by
(\ref{main1}), $\hw\hat{t}_{1,r}\hw^{-1}=\th_{i,r}$ for $r\in
J_\nu$. Thus from Lemma \ref{conj-cen-1} we have
$$
[\hat{t}_{1,r},\hat{t}_{1,s}]=\hw[\hat{t}_{1,r},\hat{t}_{1,s}]\hw^{-1}=
[\hw\hat{t}_{1,r}\hw^{-1},\hw\hat{t}_{1,s}\hw^{-1}]=[\hat{t}_{i,r},\hat{t}_{i,s}].
$$
This gives (i). Next let $J\in\supp(S)$. Then we have again from
Lemma \ref{conj-cen-1} and (\ref{main1}) that
$$
\hz_{_J}=\hw\hz_{_J}\hw^{-1}=\hw\hw_{\a_1}\hw_{\a_1+\tau_{_J}}\prod_{r\in
    J}\th_{1,r}\hw^{-1}= \hw_{\a_i}\hw_{\a_i+\tau_{_J}}\prod_{r\in
    J}\th_{i,r}.
    $$
Recall from Section \ref{EAWG} that if $\ell>1$, then $S$ is a
lattice. In particular $J\in\supp(S)$ for all $J\sub J_\nu$. This
gives (ii) and completes the proof.\qed

\begin{lem}\label{conj-cen-2}
Let $i\in J_\ell$ and $\{r,s\}\in\supp(S)$, $r<s$. Then
\begin{eqnarray*}
\hz_{_{\{r,s\}}}=\hw_{\a_i}\hw_{\a_i+\sg_r+\sg_s}\th_{i,r}\th_{i,s}=\hw_{\a_i+\sg_r+\sg_s}
\hw_{\a_i+\sg_r}\hw_{\a_i}\hw_{\a_i+\sg_s}
\end{eqnarray*}
and
\begin{eqnarray*}
\hz_{_{\{r,s\}}}^{-1}=\hw_{\a_i}\hw_{\a_i+\sg_s+\sg_r}\th_{i,s}\th_{i,r}&=&
\hw_{\a_i+\sg_s+\sg_r}\hw_{\a_i+\sg_s}\hw_{\a_i}
\hw_{\a_i+\sg_r}\\
&=&\hw_{\a_i-\sg_r+\sg_s} \hw_{\a_i-\sg_r}\hw_{\a_i}
\hw_{\a_i+\sg_s}.
\end{eqnarray*}
In particular, $\hz_{_{\{r,s\}}}=\hz_{_{\{r,s\}}}^{-1}
[\th_{i,r},\th_{i,s}]$.
\end{lem}

\proof The first equality in the statement holds by Lemma
\ref{well-define}(ii). Since by Lemma \ref{conj-cen-1}(ii),
$\hw_{\a_i}\hz_{_{\{r,s\}}}\hw^{-1}_{\a_i}=\hz_{_{\{r,s\}}}$, the
second equality also holds. Next, since $\hz_{_{\{r,s\}}}^{-1}\in
Z(\hat{\w})$ and the relations of the form (\bi) hold in
$\hat{\w}$, we have
\begin{eqnarray*}
\hz_{_{\{r,s\}}}^{-1}&=&\hw_{\a_i}\hw_{\a_i+\sg_r+\sg_s}
\hz_{_{\{r,s\}}}^{-1}\hw_{\a_i+\sg_r+\sg_s}\hw_{\a_i}\\
&=&\hw_{\a_i}\hw_{\a_i+\sg_s+\sg_r}\th_{i,s}\th_{i,r}\\
&=&\hw_{\a_i}\hw_{\a_i}\hw_{\a_i+\sg_s+\sg_r}\th_{i,s}\th_{i,r}\hw_{\a_i}\\
&=&\hw_{\a_i+\sg_r+\sg_s}
\hw_{\a_i+\sg_s}\hw_{\a_i}\hw_{\a_i+\sg_r}\\
&=&(\hw_{\a_i+\sg_s}\hw_{\a_i+\sg_r+\sg_s}
\hw_{\a_i+\sg_s})(\hw_{\a_i}\hw_{\a_i+\sg_r}\hw_{\a_i})\hw_{\a_i}\hw_{\a_i+\sg_s}\\
&=&\hw_{\a_i-\sg_r+\sg_s}\hw_{\a_i-\sg_r}\hw_{\a_i}\hw_{\a_i+\sg_s}
 \end{eqnarray*}
 This completes the proof.\qed

\begin{lem}\label{conj-11-1}
Let $(i,r),(j,s)\in J_\ell\times J_\nu$ with $r\leq s$. Then
\begin{eqnarray}\label{condition}
[\hat{t}_{i,r},\hat{t}_{j,s}]=\left\{\begin{array}{ll}
  \hz_{_{\{r,s\}}}^{(\a_i,\a_j)} & \hbox{if}\; \{r,s\}\in\supp(S),\vspace{2mm}\\
  \hz_{_{\{r,s\}}} & \hbox{if}\; \{r,s\}\not\in\supp(S).
\end{array}\right.
\end{eqnarray}
\end{lem}

\proof First, let $\{r,s\}\not\in\supp(S)$. This can happen only
if  $S$ is not a lattice. Thus    $\ell=1$ and so $i=j=1$ and the
result holds by the way $\hz_{_{\{r,s\}}}$ is defined. Next, let
$\{r,s\}\in\supp(S)$. If $i=j$, then by Lemma \ref{conj-cen-2} we
have
\begin{eqnarray*}
\hz_{_{\{r,s\}}}^2=\hz_{_{\{r,s\}}}\hz^{-1}_{_{\{r,s\}}}[\hat{t}_{i,r},\hat{t}_{i,s}]=
[\hat{t}_{i,r},\hat{t}_{i,s}].
\end{eqnarray*}
If $i\neq j$, then  using  relations (\bi) and (\bii) and Lemma
\ref{conj-cen-2} we get
 \begin{eqnarray*}
[\hat{t}_{i,r},\hat{t}_{j,s}]&=&\hat{t}_{i,r}^{-1}\hat{t}_{j,s}^{-1}\hat{t}_{i,r}\hat{t}_{j,s}\\
&=&\hw_{\a_i}\hw_{\a_i+\sg_r}\hw_{\a_j}\hw_{\a_j+\sg_s}
\hw_{\a_i+\sg_r}\hw_{\a_i}\hw_{\a_j+\sg_s}\hw_{\a_j}\\
&=& \hw_{-\a_i+\sg_r}\hw_{\a_i}\hw_{-\a_j+\sg_s}
(\hw_{\a_i}\hw_{\a_i})\hw_{\a_j}\hw_{\a_i}\hw_{-\a_i+\sg_r}
\hw_{\a_j+\sg_s}\hw_{\a_j}\\
&=&\hw_{-\a_i+\sg_r}\hw_{-w_{\a_i}(\a_j)+\sg_s}
\hw_{w_{\a_i}(\a_j)}\hw_{-\a_i+\sg_r}
\hw_{\a_j+\sg_s}\hw_{\a_j}\\
&=&(\hw_{-\a_i+\sg_r}\hw_{-w_{\a_i}(\a_j)+\sg_s}\hw_{-\a_i+\sg_r})(\hw_{-\a_i+\sg_r}
\hw_{w_{\a_i}(\a_j)}\hw_{-\a_i+\sg_r})
\hw_{\a_j+\sg_s}\hw_{\a_j}\\
&=&\hw_{_{w_{\a_i-\sg_r}(w_{\a_i}(-\a_j)+\sg_s)}}\hw_{_{w_{\a_i-\sg_r}(w_{\a_i}(\a_j))}}
\hw_{\a_j+\sg_s}\hw_{\a_j}\\
&=&\hw_{-\a_j+(\a_j,\a_i)\sg_r+\sg_s} \hw_{\a_j-(\a_j,\a_i)\sg_r}
\hw_{\a_j+\sg_s}\hw_{\a_j}\\
&=&\hw_{\a_j}\hw_{\a_j+(\a_j,\a_i)\sg_r+\sg_s}
\hw_{\a_j}\hw_{\a_j-(\a_j,\a_i)\sg_r}\hw_{\a_j}\hw_{\a_j}
\hw_{\a_j+\sg_s}\hw_{\a_j}\\
&=&\hw_{\a_j}\hw_{\a_j+(\a_j,\a_i)\sg_r+\sg_s}
\hw_{\a_j+(\a_j,\a_i)\sg_r}\hw_{\a_j} \hw_{\a_j+\sg_s}\hw_{\a_j}\\
&=&\hw_{\a_j+(\a_j,\a_i)\sg_r+\sg_s}
\hw_{\a_j+(\a_j,\a_i)\sg_r}\hw_{\a_j} \hw_{\a_j+\sg_s}.
\end{eqnarray*}
But if $(\a_i,\a_j)=-1$, then by Lemma \ref{conj-cen-2} the
expression appearing in the right of the last equality is
$\hz^{-1}_{_{\{r,s\}}}$ and is clearly $1$, if $(\a_i,\a_j)=0$.
This completes the proof.\qed

\begin{lem}\label{conj-normal}
$\hw_{\a_i}\th_{j,r}\hw_{\a_{i}}=\th_{j,r}\th_{i,r}^{-(\a_i,\a_j)}$,
$i,j\in J_\ell$, $r\in J_\nu$.
\end{lem}

\proof First, let $i=j$. Then from (\bi) and the fact that
$(\a_i,\a_i)=2$ we have
$$
\hw_{\a_i}\th_{i,r}\hw_{\a_{i}}=\hw_{\a_i}\hw_{\a_i+\sg_r}\hw_{\a_i}\hw_{\a_i}=\hw_{\a_i}\hw_{\a_i+\sg_r}
=\th_{i,r}^{-1}=\th_{i,r}\th_{i,r}^{-(\a_i,\a_i)}.
 $$
 Next let $i\neq j$. If $(\a_i,\a_j)=0$, then using (\bi) and (\bii) we
 have
 \begin{eqnarray*}
\hw_{\a_i}\th_{j,r}\hw_{\a_{i}}
&=&\hw_{\a_i}\hw_{\a_j+\sg_r}\hw_{\a_j}\hw_{\a_i}\\
&=&\hw_{w_{\a_i}(\a_j)+\sg_r}\hw_{w_{\a_i}(\a_j)}\\
&=&\hw_{\a_j+\sg_r}\hw_{\a_j}
=\th_{j,r}=\th_{j,r}\th_{i,r}^{-(\a_i,\a_i)}.
  \end{eqnarray*}
Finally, if $(\a_i,\a_j)=-1$, then using (\bi) and (\bii), we
 have
  \begin{eqnarray*}
\hw_{\a_i}\th_{j,r}\hw_{\a_i}&=&
\hw_{\a_i}\hw_{\a_j+\sg_r}\hw_{\a_j}\hw_{\a_i}\th_{i,r}^{-1}\th_{j,r}^{-1}\th_{j,r}\th_{i,r}\\
&=&\hw_{\a_i}\hw_{\a_j+\sg_r}\hw_{\a_j}\hw_{\a_i+\sg_r}\hw_{\a_j}\hw_{\a_j+\sg_r}\th_{j,r}\th_{i,r}\\
&=&\hw_{\a_i}\hw_{\a_j+\sg_r}\hw_{\a_i+\a_j+\sg_r}\hw_{\a_j+\sg_r}\th_{j,r}\th_{i,r}\\
&=&\hw_{\a_i}\hw_{\a_i}\th_{j,r}\th_{i,r}=\th_{j,r}\th_{i,r}=\th_{j,r}\th_{i,r}^{-(\a_i,\a_i)}.
 \end{eqnarray*}
This completes the proof.\qed

 \begin{lem}\label{cofi-two}
$\hw_{\a_1+2\sum_{r\in J}\sg_r}=
 (\prod_{r\in J} \th_{1,r})\hw_{\a_1}(\prod_{r\in J}
\th_{1,r})^{-1}$,\quad $J\subseteq J_\nu$.
\end{lem}

\proof Let $w=\prod_{r\in J} t_{1,r}$. From (\ref{prod1}) we have
$w(\a_1)=\a_1+\sum_{r\in J}2\sg_r$, and so by (\bii) and
(\ref{main1}) we get
\begin{eqnarray*} \qquad\qquad\qquad\hw_{\a_1+2\sum_{r\in
J}\sg_r}=\hw_{w(\a_1)}=(\prod_{r\in J}
\th_{1,r})\hw_{\a_1}(\prod_{r\in J}
\th_{1,r})^{-1}.\mbox{\qquad\qquad\qquad\qed}
\end{eqnarray*}

\begin{lem}\label{spilit-con-1}
 $\hz_{_J}^2\in \la \hz_{_{\{r,s\}}}\mid r,s\in J_\nu,\;r<s\ra$, $J\subseteq J_\nu$.
\end{lem}
\proof If $J\in\supp(S)$ and $\tau_{_J}=\sum_{r\in J}\sg_r$, then
using relations (\bi), (\bii) and Lemmas \ref{cofi-two} and
\ref{conj-normal} we get
\begin{eqnarray*}
\hz_{_J}^2&=&\hz_{_J}\hz_{_J}=(\hw_{\a_1}\hw_{\a_1+\tau_{_J}}\prod_{r\in
J}\th_{1,r})\hz_{_J}\\
&=&\hw_{\a_1}\hw_{\a_1+\tau_{_J}}\hz_{_J} \prod_{r\in
J}\th_{1,r}\\
&=&\hw_{\a_1}(\hw_{\a_1+\tau_{_J}}\hw_{\a_1}\hw_{\a_1+\tau_{_J}})
\prod_{r\in J}\th_{1,r}\prod_{r\in J}\th_{1,r}\\
&=&\hw_{\a_1}\hw_{\a_1+2\tau_{_J}}
\big(\prod_{r\in J}\th_{1,r}\big)^2\\
&=&\hw_{\a_1}(\prod_{r\in J} \th_{1,r})\hw_{\a_1}(\prod_{r\in J}
\th_{1,r})^{-1}
\big(\prod_{r\in J}\th_{1,r}\big)^2\\
&=&\prod_{r\in J} \th_{1,r}^{-1}\prod_{r\in J} \th_{1,r}.
\end{eqnarray*}
But this latter belongs to $\la\hz_{_{\{r,s\}}}\mid r,s\in
J_\nu,\;r<s\ra$, by Lemma \ref{conj-11-1}. If $J\not\in supp(S)$,
then the result is clear by the way $z_{_J}$ is defined and Lemma
\ref{conj-11-1}.\qed
\begin{lem}\label{free}
$\la\hz_{_{\{r,s\}}}\mid 1\leq r<s\leq\nu\ra$ is a free abelian
group of rank $\nu(\nu-1)/2$.
 \end{lem}
\proof Let $\prod_{1\leq
r<s\leq\nu}\hz_{_{\{r,s\}}}^{m_{r,s}}=1$,\; $m_{r,s}\in\bbbz$. By
(\ref{psi}) we have that
$$1=\psi(\prod_{1\leq
r<s\leq\nu}\hz_{_{\{r,s\}}}^{m_{r,s}})=\prod_{1\leq
r<s\leq\nu}z_{_{\{r,s\}}}^{m_{r,s}}.$$ Thus from Lemma
\ref{indepen}, it follows that $m_{r,s}=0$, for any $1\leq
r<s\leq\nu$.\qed

\begin{lem}\label{tavan-2-z}
If $J\in\supp(S)$, then $$\hz_{_J}^2=\prod_{\{r,s\in
J~|~r<s\}}\hz_{_{\{r,s\}}}^{3-\d(r,s)}.$$
\end{lem}
\proof By (\ref{psi}), Lemmas \ref{free} and \ref{indepen}, the
restriction $\psi$ to $\la \hz_{_{\{r,s\}}}~|~1\leq r<s\leq\nu\ra$
induces the isomorphism
\begin{eqnarray}\label{BBB}
\la \hz_{_{\{r,s\}}}~|~1\leq
r<s\leq\nu\ra\cong^{^{\hspace{-2mm}\psi}}\la
z_{_{\{r,s\}}}~|~1\leq r<s\leq\nu\ra.
\end{eqnarray}
From  the way $z_{_J}$ is defined and  (\ref{final1}), it follows
that
$$z_{_J}^2=\prod_{\{r,s\in J~|~r<s\}}z_{_{\{r,s\}}}^{3-\d(r,s)}.$$
Then from Lemma \ref{spilit-con-1}, the facts (\ref{psi}) and
(\ref{BBB}) we get
\begin{eqnarray*}
\hz_{_J}^2&=&\psi^{-1}(z_{_J}^2)\\
&=&\psi^{-1}(\prod_{\{r,s\in
J~|~r<s\}}z_{_{\{r,s\}}}^{3-\d(r,s)})\\
&=&\prod_{\{r,s\in
J~|~r<s\}}\psi^{-1}(z_{_{\{r,s\}}}^{3-\d(r,s)})=\prod_{\{r,s\in
J~|~r<s\}}\hz_{_{\{r,s\}}}^{3-\d(r,s)}.
\end{eqnarray*}
\qed
\begin{lem}\label{form>1}
If $\ell>1$, then each element $\hw$ of $\hat{\w}$ has a unique
expression in the form
 \begin{eqnarray}\label{A>1}
\hw=\hw(\hat{\dot{w}},n_{i,r},m_{r,s}):=\hat{\dot{w}}
\prod_{r=1}^{\nu}\prod_{i=1}^{\ell}\th_{i,r}^{n_{i,r}}\prod_{1\leq
r<s\leq\nu}\hz_{_{\{r,s\}}}^{m_{r,s}},
\end{eqnarray}
where $n_{i,r},\;m_{r,s}\in\bbbz$ and
$\hat{\dot{w}}\in\hat{\dot{\w}}$.
 \end{lem}

\proof Let $\hat{w}\in\hat{\w}$.  By Lemmas
\ref{generators-con-1}, \ref{conj-cen-1}, \ref{conj-11-1}, and
\ref{conj-normal}, $\hat{w}$ can be written in the form
(\ref{A>1}).  Let $\hw(\hat{\dot{w}}',n_{i,r}',m_{r,s}')$
 be another expression of
$\hw$ in the form (\ref{A>1}). Then
$\psi(\hw(\hat{\dot{w}},n_{i,r},m_{r,s}))=\psi(\hw(\hat{\dot{w}}',n_{i,r}',m_{r,s}')$
and so from (\ref{psi}), (\ref{finite-case}) and Proposition
\ref{unique-6}, we get $\hat{\dot{w}}=\hat{\dot{w}}'$ and
$n_{i,r}=n_{i,r}'$, $(i,r)\in J_\ell\times J_\nu$ and
$m_{r,s}=m_{r,s}'$ for all $r,s\in J_\nu$. \qed

We now are able to give a new proof to the following theorem which
is due to \cite{K}. Indeed, the proof given in \cite{K} is based
on the fact that the presentation by conjugation holds for finite
and affine Weyl groups while in our proof we have only used this
fact for the finite Weyl group. In particular, the proof given
here can also apply to the affine case.

\begin{thm}\label{form-conj-1a}
Let $R$ be a simply laced extended affine root system of rank
$\ell>1$. Then the Weyl group $\w$ of $R$ has the presentation by
conjugation.
\end{thm}

\proof  It is enough to show that the epimorphism $\psi$ defined
by (\ref{main}) is one to one. Let $\hat w\in\hat{\w}$ so that
$\psi(\hat w)=1$. By Lemma \ref{form>1}, $\hat w$ has an
expression of the form (\ref{A>1}).  Then by (\ref{psi}) we have
that
\begin{eqnarray*}
1=\psi(\hw)=\psi(\hat{\dot{w}})
\prod_{r=1}^{\nu}\prod_{i=1}^{\ell}t_{i,r}^{n_{i,r}}\prod_{1\leq
r<s\leq\nu}z_{_{\{r,s\}}}^{m_{r,s}}.
 \end{eqnarray*}
Therefore from (\ref{finite-case}) and Proposition \ref{unique-6},
it follows that $\hat{\dot{w}}=1$ and $n_{i,r}=0$, $(i,r)\in
J_\ell\times J_\nu$ and $m_{r,s}=0$ for all $r,s\in J_\nu$ and so
$\hat w=1$. \qed

\section{\bf PRESENTATION BY CONJUGATION FOR TYPE $A_1$}\label{type-1}
\setcounter{equation}{0}

This section contains the main results of our work. In particular,
we give the necessary and sufficient conditions for an extended
affine Weyl group of type $A_1$ to have the presentation by
conjugation (Theorem \ref{reduced}). Throughout this section $R$
is of type $A_1$, that is
$$R=(S+S)\cup(\pm\a_1+S).$$ As in the previous sections let $\w$ be the
extended affine Weyl group of $R$ and $\hat{\w}$ be the group
defined in Definition \ref{DEF}. We recall from (\ref{pi}) and
Lemma \ref{base-1} that the set $$\Pi=\{\a_1+\tau_{_J}\mid
J\in\supp(S)\}$$ satisfies $\w=\la w_\a\mid\a\in\Pi\ra$ and $
\w\Pi=\rtimes$.

\begin{lem}\label{minimal}
 $\{\hw_{\a}\mid\a\in\Pi\}$
 is a minimal set of generators for $\hat{\w}$.
 \end{lem}

 \proof By Lemma \ref{generators-con-1}(i),
 the set in the statement generates $\hat{\w}$. To show that
 it is minimal fix $\b=\a_1+\tau_{_{J_0}}\in\Pi$,  $J_0\in\supp(S)$. We have
$\hat{\w}=F/N$, where $F$ is the free group on the set
$\{r_\a\mid\a\in \r^\times\}$ and
 $N$ is the normal closure of the set $\{r_\a^2,\; r_\a r_\b r_\a r_{w(\b)}\mid\a,\b\in\r^\times\}$
 in $F$. Then $\hw_{\a}=r_\a N$, $\a\in\rtimes$. From  (\ref{union}), it follows
 that
 \begin{eqnarray*}
\r^\times=\biguplus_{J\in\supp(S)}\r_{_{J}},\quad\mbox{where}\quad
\r_{_{J}}:=\pm\a_1+\tau_{_{J}}+2\Lam
 \end{eqnarray*}
 Since $(\a,\a)=2$ for $\a\in\rtimes$, we have
 \begin{eqnarray}\label{stabel}
\w\r_{_{J}}\subseteq \r_{_{J}}\qquad(J\in\supp(S)).
 \end{eqnarray}
Now let $\varphi:F\longrightarrow\bbbz_2$ be the epimorphism
induced by the assignment
 $$\varphi(r_{\a})=\left\{\begin{array}{ll}
   1, & \hbox{if}\;\; \a\in \r_{_{J_0}},\vspace{2mm} \\
   0, & \hbox{if}\;\; \a\in\r^\times\setminus\r_{_{J_0}}.
 \end{array}\right.$$
From (\ref{stabel}), it follows that $\varphi(r_\a^2)=0$ and
$\varphi(r_\a r_\b r_\a
  r_{w_\a(\b)})=0$,  for any $\a,\b\in\r^\times$  and  so
  $\varphi(N)=\{0\}$. Thus the epimorphism $\varphi$ induces a
  unique epimorphism
  $\bar{\varphi}:\hat{\w}\longrightarrow\bbbz_2$ so that
  $\bar{\varphi}(\hw_\a)=\varphi(r_\a)$, $\a\in\r^\times$. So $\hw_\b\not=\hw$ for any
  $\hw\in\la\hw_\a\mid\a\in\Pi\setminus\b\ra$, as
 $\bar{\psi}(\hw_\b)=1$ and $\bar{\psi}(\hw)=0$.\qed

\begin{lem}\label{uique} If
$\prod_{1\leq
r<s\leq\nu}\hz_{_{\{r,s\}}}^{m_{r,s}}\prod_{J\in\esupp(S)}
\hz_{_{J}}^{\epsilon_{_J}}=1, $ where all $m_{r,s}$'s are integers
and $\ep_{_J}\in\{0,\pm 1\}$ for all $J$. Then
$m_{r,s}=0=\ep_{_J}$ for all $r,s,J$.
\end{lem}

\proof  By Lemma \ref{free}, it is enough to show that
$\ep_{_J}=0$ for all $J\in\esupp(S)$. Suppose to the contrary that
$\ep_{_{J_0}}\neq0$ for some $J_0\in\esupp(S)$. Then
 $$\hz_{_{J_0}}^{-\ep_{_{J_0}}}=\prod_{1\leq
r<s\leq\nu}\hz_{_{\{r,s\}}}^{m_{r,s}}\prod_{J\in\esupp(S)\setminus\{J_0\}}
\hz_{_{J}}^{\epsilon_{_J}},$$ and so $\hz_{_{J_0}}\in\la
\hz_{_{\{r,s\}}},\hz_{_{J}}\mid 1\leq
r<s\leq\nu,\;J\in\esupp(S)\setminus\{J_0\}\ra\ $. Therefore from
the way $\hz_{_J}$'s  are defined (see (\ref{test-1})) it follows
that $\hw_{\a_1+\tau_{_{J_0}}}\in\la \hw_{\a_1+\tau_{_J}}\mid
J\in\supp(S)\setminus\{J_0\}\ra$ and this contradicts Lemma
\ref{minimal}. \qed

From now on we consider pairs of collections of the form
\begin{equation}\label{form0}
(\mb,\epb)=\big(\{m_{r,s}\}_{1\leq
r<s\leq\nu},\{\ep_{_J}\}_{J\in\esupp(S)}\big), \quad
m_{r,s}\in\bbbz,\;\;\ep_{_J}\in\{0,1\},
\end{equation}
where as before if $\esupp(S)=\emptyset$, we interpret $\epb$ as
the zero collection. To each such a pair we assign a central
element of $\hat{\w}$ by
\begin{equation}\label{reduced-form}
\hu(\mb,\epb)=\prod_{1\leq
r<s\leq\nu}\hz_{_{\{r,s\}}}^{m_{r,s}}\prod_{J\in\esupp(S)}
\hz_{_{J}}^{\epsilon_{_J}}.
\end{equation}
From Lemma \ref{uique} we have
\begin{equation}\label{unique1}
\hu(\mb,\epb)=1\Longleftrightarrow (\mb,\epb)=(\{0\},\{0\}).
\end{equation}

\begin{lem}\label{form-conj-1}
Each element $\hw$ of $\hat{\w}$ has a unique expression of the
form
 \begin{eqnarray}\label{A-1}
\hw=\hw(n,n_{r},\mb,\epb):=\hw_{\a_1}^n
\prod_{r=1}^{\nu}\th_{1,r}^{n_{r}}\;\hu (\mb,\epb),
\end{eqnarray}
where $n\in\{0,1\}$, $n_{r}$'s are integers and $(\mb,\epb)$ is of
the form (\ref{form0}).
 \end{lem}

\proof First we can express each element $\hw\in\hat{\w}$ in terms
of the generators given in Lemma \ref{generators-con-1}(ii). Next
we can reorder the appearance of generators in any such expression
using Lemmas \ref{conj-cen-1}, \ref{conj-11-1}, and
\ref{conj-normal}. Finally we get the appropriate expression using
Lemmas \ref{conj-cen-1}(ii), \ref{conj-11-1} and
\ref{spilit-con-1}. Now to complete the proof it is enough to show
that the expression of $\hw$ in the form (\ref{A-1}) is unique.
Let $\hw(n',n_{r}',\mb',\epb')$
 be another expression of
$\hw$ in the form (\ref{A-1}) where $\mb'=\{m'_{r,s}\}$ and
$\epb'=\{\ep'_{_J}\}$. Then applying the homomorphism $\psi$ on
these two expressions of $\hw$ and using Proposition
\ref{unique-6} we get $n=n'$ and $n_{r}=n_{r}'$ for all $r\in
J_\nu$. Thus $\hu(\mb,\epb)=\hu(\mb',\epb').$ Then from Lemma
\ref{conj-cen-1}(ii), we have that
$$\prod_{1\leq
r<s\leq\nu}\hz_{_{\{r,s\}}}^{m_{r,s}-m_{r,s}'}\prod_{J\in\esupp(S)}
\hz_{_{J}}^{\epsilon_{_J}-\epsilon_{_J}'}=1,$$
 and so from Lemma \ref{uique}, we get
 that $m_{r,s}=m'_{r,s}$ and
$\epsilon_{_J}=\epsilon_{_J}'$  for all   $1\leq r<s\leq\nu$ and
$J\in\esupp(S)$. \qed

 Set \begin{eqnarray}\label{HZ}
\hat{Z}:=\la\hz_{_J}\mid J\sub J_\nu\ra\leqslant Z(\hat{\w}).
\end{eqnarray}
The following Corollary follows immediately from Lemma
\ref{form-conj-1}.

 \begin{cor}\label{form-center}
 Each element $\hu$ of $ \hat{Z}$ has a
unique expression of the form
$$\hu=\hu\big(\mb,\epb),
$$
where $(\mb,\epb)$ is of the form (\ref{form0}).
 \end{cor}

\begin{pro}\label{redu-3}
(i) $\ker(\psi)\subseteq \hat{Z}$.

(ii) If $\hu\in\hat{Z}$, then the following statements are
equivalent:

\begin{itemize}
\item[(a)] $\hu\in\ker(\psi)$,

\item[(b)] $\hu^2=1$,

\item[(c)] $|\hu|<\infty$,

\item[(d)] $\hu=\hu(\mb,\epb)$, where $\epb$ is an integral
collection satisfying
\begin{equation}\label{mrs}
m_{r,s}=-\frac{1}{\d(r,s)}\sum_{J\in\esupp(S)}
\d(_J,r,s)\epsilon_{_J},\qquad (1\leq r<s\leq\nu).
\end{equation}
\end{itemize}

(iii) The assignment $\epb=\{\ep_{_J}\}\longmapsto\hu(\mb,\epb)$,
where $\mb=\{m_{r,s}\}$ satisfies (\ref{mrs}), is a one to one
correspondence from the set of integral collections for $S$ onto
$\ker(\psi)$.
\end{pro}

\proof (i) Let $\hw\in\ker(\psi)$. By  Lemma \ref{form-conj-1},
$\hw$
 has a unique expression of the form (\ref{A-1}).
 Then applying the homomorphism $\psi$, we get
 $$
1=\psi(\hw)=\psi(\hw_{\a_1}^n)\prod_{r=1}^{\nu}\psi(\th_{1,r}^{n_{r}})
\psi\big(\hu(\mb,\epb)\big)=w_{\a_1}^n\prod_{r=1}^{\nu}t_{1,r}^{n_{r}}\;z,
$$
where $z=\psi\big(\hu(\mb,\epb)\big)\in Z(\w)$. Then from
Proposition \ref{unique-6} we obtain $n=0=n_{r}$ for all $r\in
J_\nu$ and so $\hw=\hu(\mb,\epb)\in \hat{Z}$.

(ii) Let $\hu\in\hat{Z}$. By Corollary \ref{form-center},
$\hu=\hu(\mb,\epb)$ for a pair $(\mb,\epb)$ of the form
(\ref{form0}) . Let $n\in\mathbb Z_{\geq0}$. By (\ref{psi}),
(\ref{def-z}), (\ref{free-abelian}) and  (\ref{final1}) we have
\begin{eqnarray*}
\psi(\hu^n)&=&\prod_{1\leq
r<s\leq\nu}z_{_{\{r,s\}}}^{nm_{r,s}}\prod_{J\in\esupp(S)}
z_{_{J}}^{n\epsilon_{_J}}\\
&=&\prod_{1\leq
r<s\leq\nu}c_{r,s}^{\d(r,s)nm_{r,s}}\prod_{J\in\esupp(S)}
\prod_{\{r,s\in J\mid r<s\}}c_{r,s}^{n\ep_{_J}}\\
&=& \prod_{1\leq r<s\leq\nu}c_{r,s}^{\d(r,s)nm_{r,s}} \prod_{1\leq
r<s\leq\nu}\prod_{J\in\esupp(S)}c_{r,s}^{\d(_J,r,s)n\ep_{_J}}\\
&=& \prod_{1\leq r<s\leq\nu}c_{r,s}^{\d(r,s)nm_{r,s}} \prod_{1\leq
r<s\leq\nu}c_{r,s}^{\sum_{J\in\esupp(S)}\d(_J,r,s)n\ep_{_J}}\\
&=& \prod_{1\leq
r<s\leq\nu}c_{r,s}^{n\d(r,s)m_{r,s}+n\sum_{J\in\esupp(S)}\d(_J,r,s)\ep_{_J}}.
\end{eqnarray*}
Then it is immediate from (\ref{free-abelian}) that
$(c)\Longrightarrow (d)\Longleftrightarrow (a)$. Since
$(b)\Longrightarrow (c)$, it remains to show that
$(a)\Longrightarrow (b)$. Let $\hu\in \ker(\psi)$. By Lemma
\ref{spilit-con-1}, we have
\begin{eqnarray}\label{tavan2}
\hu^2=\prod_{1\leq r<s\leq\nu}\hz_{_{\{r,s\}}}^{m_{r,s}},
\end{eqnarray}
for some $m_{r,s}\in\bbbz$. Then by (\ref{psi}), we have that
$$1=\psi(\hu^2)=\prod_{1\leq
r<s\leq\nu}z_{_{\{r,s\}}}^{m_{r,s}}.$$ Now Lemma \ref{indepen}
gives $m_{r,s}=0$ for all $r,s$ and so $\hu^2=1$.

(iii) This is an immediate consequence of Corollary
\ref{form-center} and part (ii).\qed

\begin{cor}\label{fgag}
$\ker(\psi)$ is isomorphic to a direct sum of at most
$|\esupp(S)|$-copies of $\bbbz_2$.
\end{cor}

\proof By Proposition \ref{redu-3}, each non-trivial element of
$\ker(\psi)$ if of order $2$ and has the form $\hu(\mb,\epb)$,
where $\epb$ is an integral collection satisfying (\ref{mrs}). But
there are at most $2^{|\esupp(S)|}$ integral collections (see
Section \ref{Semilattices}). Now the result is clear as
$\ker(\psi)$ is abelian.\qed

Let $n_0\in\bbbz_{\geq 0}$ be the number of copies of $\bbbz_2$
involved in $\ker(\psi)$, then
$$
|\ker(\psi)|=2^{n_0}.$$

\begin{cor}\label{lattice-1}
If $\{r,s\}\in\supp(S)$ for all $1\leq r<s\leq\nu$, then
$n_0=|\esupp(S)|$. In particular, if $S$ is a lattice then
$n_0=2^{\nu}-1-\nu-\left(
\begin{array}{c}
  \nu \\
  2 \\
\end{array}\right)$.
\end{cor}

\proof By Example \ref{exa1}(ii), under conditions in the
statement, there is exactly $2^{n_0}$ integral collections. Now
the result follows from Proposition \ref{redu-3}(iii).\qed

\begin{thm}\label{reduced}
 Let $R=R(A_1, S)$ be an extended affine
root system of type $A_1$ with extended affine Weyl group $\w$.
Then the following statements are equivalent:

 (i)  $\hat{\w}\cong\w$ ($\w$ has the presentation by conjugation).

 (ii) $Z(\hat{\w}){\cong} Z(\w)$.

(iii) $Z(\hat{\w})$ is a free abelian group.

(iv) $\psi:\hat{\w}\longrightarrow\w$ is an isomorphism.

(v) The trivial collection is the only integral collection for
$S$.

(vi) $\{w_{\a}\mid \a\in\Pi\}$ is a minimal set of generators for
$\w$, where $\Pi$ is given by (\ref{pi}).
\end{thm}

\proof Clearly we have (iv)$\Longrightarrow$(i)$\Longrightarrow$
(ii). By Proposition \ref{conj-6}, (ii) $\Longrightarrow$ (iii).
By Proposition \ref{redu-3}, (iii) implies (iv) and (iv) and (v)
are equivalent. From Lemma \ref{minimal} and (\ref{psi}) it
follows that (iv) implies (vi). So it remains to show that
(vi)$\Longrightarrow$(i). Let $\hu\in\ker(\psi)$. By Proposition
\ref{redu-3}, $\hu$ has the form  $\hu(\mb,\epb)$, where $\epb$ is
an integral collection satisfying (\ref{mrs}). Now we claim  that
$\ep_{_J}=0$ for all $J\in\esupp(S)$ and so $\hu=1$. Suppose to
the contrary that $\ep_{_{J_0}}\neq0$ for some $J_0\in\esupp(S)$.
Then from (\ref{psi}) we have
\begin{eqnarray*}
1=\psi(\hu)=\prod_{1\leq
r<s\leq\nu}z_{_{\{r,s\}}}^{m_{r,s}}\prod_{J\in\esupp(S)}
z_{_{J}}^{\epsilon_{_J}}\\
\end{eqnarray*}
and so
 $$z_{_{J_0}}^{-\ep_{_{J_0}}}=\prod_{1\leq
r<s\leq\nu}z_{_{\{r,s\}}}^{m_{r,s}}\prod_{J\in\esupp(S)\setminus\{J_0\}}
z_{_{J}}^{\epsilon_{_J}},$$ thus $z_{_{J_0}}\in\la
z_{_{\{r,s\}}},z_{_{J}}\mid 1\leq
r<s\leq\nu,\;J\in\esupp(S)\setminus\{J_0\}\ra\ $. Therefore from
Lemma \ref{conj-sim}(i) it follows that
$w_{\a_1+\tau_{_{J_0}}}\in\la w_{\a_1+\tau_{_J}}\mid
J\in\supp(S)\setminus\{J_0\}\ra$ and this contradicts (vi) and so
(vi)$\Longrightarrow$(i).\qed

\begin{cor}\label{reduced-3}
Let $R=R(A_1, S)$ and $\w$ be as in Theorem \ref{reduced} and
either of the following holds:

 (i) There exists  $J\in\esupp(S)$ such that
$\{r,s\}\in\supp(S)$ for all $r,s\in J$.

(ii) $S$ is a lattice of rank $\geq 3$.

(iii)  $\rank(S)>3$ and $\ind(S)=2^\nu-2$.

Then $\w$ does not have the presentation by conjugation.
\end{cor}

 \proof By Example \ref{exa1}, under either of conditions (i)-(iii)
 in the statement, there is at least one non-trivial integral
 collection for $S$. Now the result follows from Theorem
 \ref{reduced}.\qed

\begin{cor}\label{exist}
Let $\nu\geq3$ and  $\nu+4\leq m\leq2^\nu-1$. Then there exists an
extended affine root system $R=R(A_1,S)$ of type $A_1$ and nullity
$\nu$ with $\ind(S)=m$ such that its extended affine Weyl group
dose not have the presentation by conjugation.
\end{cor}

\proof Let $J^1,J^2,\ldots,J^m$ be  distinct subsets of $J_\nu$
such that $J^i=\{i\}$, $1\leq i\leq \nu$, $J^{\nu+1}=\{1,2\}$,
$J^{\nu+2}=\{1,3\}$ , $J^{\nu+3}=\{2,3\}$ and
$J^{\nu+4}=\{1,2,3\}$. Set
$$\Lam:=\bbbz\sg_1\oplus\cdots\oplus\bbbz\sg_\nu \andd
S:=\bigcup_{J\in\{\emptyset,J^1,\ldots,J^m\}}(\tau_{_J}+2\Lam).$$
Then $S$ is a semilattice of rank $\nu$ with $\ind(S)=m$ and
$\supp(S)=\{\emptyset,J_1,\ldots,J_m\}$. Let $\rd=\{\pm\a_1\}$ be
a finite root system of type $A_1$. Then $R=(S+S)\cup(\rd+S)$  is
an extended affine root system of type $A_1$ of nullity $\nu$ and
index $m$. By Corollary \ref{reduced-3}(i), the extended affine
Weyl group of $R$ dose not have the presentation by
conjugation.\qed

\begin{cor}\label{main-22}
Let $R=R(A_1, S)$ and $\w$ be as in Theorem \ref{reduced}, $\nu$
the nullity of $R$ and $m=\ind(S)$. If either of the following
statements holds, then $\w$ has the presentation by conjugation:

(i) $m\leq\nu+3$,

(ii) $\nu\leq3$ and $m\neq7$,

(iii) $\esupp(S)=\emptyset$.

 In particular, if  $\nu\leq 3$, then $\w$ has the
 presentation by conjugation if and only if
 $m\not=7$.
\end{cor}

\proof By Theorem \ref{reduced}, we only need to show that the
only integral collection for $S$ is the trivial collection. To see
this, use Lemma \ref{reduced-1} for part (i). Part (ii) is a
consequence of (i) as $m\leq 7$ when $\nu\leq 3$. For part (iii)
see Example \ref{exa1}. Next, note that if $\nu\leq 2$, then
$\esupp(S)=\emptyset$ and so the last statement holds by part
(iii). Finally, if $\nu=3$ then the result holds by part (ii) and
Corollary \ref{exist}.\qed

\begin{cor}\label{thm-2}
Let $R=(S+S)\cup(\rd+S)$ and $R'=(S'+S')\cup(\rd+S')$ be two
extended affine root systems of type $A_1$, where $S$ and $S'$
have the same $\bbbz$-span and $S\sub S'$ . Then if the extended
affine Weyl group of $R'$ has the presentation by conjugation,
then also has the extended affine Weyl group of $R$.
\end{cor}

\proof It follows immediately from Theorem \ref{reduced} and
Example \ref{exa1}(iv).\qed

\begin{rem}\label{rem1}
(i) In \cite[Remark 2.14]{A3} the author, by suggesting an
example, has conjectured that certain EAWG's may not have the
presentation by conjugation. The Corollaries \ref{reduced-3} and
\ref{exist}, give an affirmative answer to this conjecture.
However, as it is shown in Corollary \ref{main-22}, the suggested
example in \cite{A3} is not an appropriate example. In fact, the
root system in the mentioned example is an $A_1$-type extended
affine root system of nullity $3$ and index $6$, and so by
Corollary \ref{main-22}, its Weyl group has the presentation by
conjugation.

(ii) One can easily see that Corollary \ref{thm-2} is false if we
remove the assumption  $\la S\ra=\la S'\ra.$
\end{rem}

 \end{document}